\documentclass[3p]{elsarticle}

\journal{N.N.}

\usepackage{amssymb,latexsym,amsmath,epsfig,amsthm} 

\usepackage[linesnumbered,ruled,noline]{algorithm2e}
\SetKwIF{If}{ElseIf}{Else}{if}{}{else if}{else}{end}%
\SetKwFor{For}{for}{}{end}
\SetKwFor{While}{while}{}{end}
\SetKwProg{Fn}{function}{}{end}
\SetKwComment{Comment}{ /\!/\,}{}
\SetCommentSty{itshape}

\newcommand{\coeff}{\mathrm{coeff}}
\newcommand{\ee}{\mathrm{e}}
\renewcommand{\AA}{\mathtt{A}}
\newcommand{\BB}{\mathtt{B}}
\newcommand{\QQ}{\mathbb{Q}}

\newtheorem{theorem}{Theorem}

\begin{document}

\begin{frontmatter}
    \title{A simple and efficient algorithm for computing the Baker--Campbell--Hausdorff series}
    
    \author{Harald Hofst\"atter}
    \address{Reitschachersiedlung 4/6, 7100 Neusiedl am See, Austria}
    \ead{hofi@harald-hofstaetter.at}
    \ead[url]{www.harald-hofstaetter.at}


    \begin{abstract}A new algorithm for computing coefficients of the Baker--Campbell--Hausdorff series is presented,         
        which can be straightforwardly implemented in any general-purpose programming language
        or computer algebra system. The algorithm avoids
        rational arithmetic and performs all its calculations in pure integer arithmetic, 
        allowing for a very efficient implementation.   
        An implementation in the Julia programming language is available. 
    \end{abstract}

\end{frontmatter}

\section{Introduction}
The Baker--Campbell--Hausdorff (BCH) series is defined as the element 
\begin{equation*}
H = \log(\ee^{\AA}\ee^{\BB})
= \sum_{k=1}^\infty\frac{(-1)^{k+1}}{k}\big(\ee^{\AA}\ee^{\BB}-1\big)^k
= \sum_{k=1}^\infty\frac{(-1)^{k+1}}{k}\bigg(\sum_{i+j>0}\frac{1}{i!j!}\AA^i\BB^j\bigg)^k 
\end{equation*}
in the ring $\QQ\langle\langle\AA,\BB\rangle\rangle$  of formal power series 
in the non-commuting variables $\AA$ and $\BB$ with rational coefficients.
The BCH series can be written as a sum $H=\sum_{n=1}^{\infty}H_n$
of homogeneous components 
\begin{equation}\label{eq:Hn}
H_n= \sum_{w\in\{\AA,\BB\}^n} h_w w,\quad n=1,2,\dots,
\end{equation}
where $\{\AA,\BB\}^n$
denotes the finite set of all words $w=w_1\cdots w_n$ ($w_i\in\{\AA,\BB\}$) of length (degree) $n$ over the 
alphabet $\{\AA, \BB\}$, and
$h_w$ denotes the coefficient of such a word in the BCH series $H$.

A classical result known as the
Baker--Campbell--Hausdorff theorem (see e.g.~\cite{HHproof})
 states that the homogeneous components $H_n$ are  \emph{Lie polynomials} 
 which means that they
can be written as linear combinations of 
$\AA$ and $\BB$ and (possibly nested) 
 commutator terms in $\AA$ and $\BB$. 
 An explicit representation of $H_n$ as a Lie polynomial is given by
\begin{equation}\label{eq:dynkin}
    H_n= \frac{1}{n}\sum_{w\in\{\AA,\BB\}^n} h_w [w], 
\end{equation}
with the same coefficients $h_w$ as in (\ref{eq:Hn}),
where $[w]$ denotes the iterated commutator
\begin{equation*}
    [w] = [w_1,[w_2,[\dots,[w_{n-1},w_n]\dots]]]
\end{equation*}
built from the word $w=w_1\cdots w_n\in\{\AA,\BB\}^n$, see \cite{Thompson82}.
It should be stressed that for $n\geq 2$, (\ref{eq:dynkin}) is not a
representation of $H_n$ as a linear combination of {\em linearly independent} 
commutators. 
Such more compact representations can e.g.~be calculated with the algorithm used
in \cite{HHbch}, which needs the computation of 
(a subset of) the  coefficients $h_w$ as a preliminary step.

Several methods for computing the coefficients $h_w$ have been proposed 
in the  literature.
Some of their implementations require only very few lines of code in a computer
algebra system.
For example, Weyrauch and Scholz \cite{WeyrauchScholz} provide a very concise  Mathematica implementation
of a method due to Goldberg \cite{G}, and 
Van--Brunt and Visser \cite{VanBruntVisser} provide a very concise Maple 
implementation of a method which is based on the algorithm of Reinsch \cite{Reinsch}.
It should be emphasized that these implementations make significant use of 
the power of the computer algebra system, thereby hiding the complexity 
of the respective methods.

In this paper, we present a new algorithm for the computation of the BCH coefficients $h_w$, which is
based on only the most basic programming constructs, and which uses only basic integer arithmetic.
The complexity of the algorithm is thus explicitly visible, but, on the  other hand, 
this allows for a straightforward and efficient implementation in any reasonable programming
language or computer algebra system.
An implementation in the Julia programming language is available at 
\cite{HHjulia}. We discuss the efficiency of this implementation in 
Section~\ref{Sect:performance}.
Note that depending on the available integer data types, there may be a limit on 
the  degrees $n$  of the coefficients $h_w$ (i.e., on the lengths of the words 
$w$) that our algorithm can calculate. 
For standard 64-bit integers, which are available for most current compilers for the Fortran and
the C programming languages, the limit is $n\leq 19$, and for
128-bit integers, which e.g.~are available 
for many compilers for the C programming language on modern computer architectures,
this limit is $n\leq 30$.
Higher degrees may require a library for multi-precision integer arithmetic.

    %



\section{Algorithm for computing coefficients of the BCH series}
Although the BCH coefficients are rational numbers, 
our algorithm performs all its calculations for a coefficient $h_w$ 
in pure integer arithmetic, except for a final division of the calculated
numerator of $h_w$ by an a priori known denominator.
For this to be possible  some a-priori information about the numerators 
of the coefficients is necessary,
which is provided by the following theorem proved in \cite{HHdenom, HHsmallest}.
\begin{theorem}\label{Thm:MainTheorem}
    For $n\geq 1$ define
    \begin{equation}\label{eq:d_n}
    d_n = \prod_{p\ \mathrm{prime}, \ p<n}p^{\max\{t:\  p^t\leq s_p(n)\}}, 
    \end{equation}
    where $s_p(n)=\alpha_0+\alpha_1+\ldots+\alpha_r$ denotes  the sum of
    the digits in the $p$-adic expansion  $n=\alpha_0+\alpha_1p+\ldots+\alpha_rp^r$.
    Then $n!\,d_n$ is the smallest common denominator for all coefficients of words
    of length $n$ in the Baker--Campbell--Hausdorff series $H=\log(\ee^{\AA}\ee^{\BB})$.
    \end{theorem}
\begin{table}\label{table:d_n}
\caption{Values of $d_n$ such that $n!d_n$ is the smallest common denominator of all
coefficients of degree $n$ in the BCH series.}
\begin{center}
\begin{tabular}{c|rrrrrrrrrrrrrrr}
\hline
$n$ &
1& 2& 3& 4& 5& 6& 7& 8&  9& 10& 11& 12&  13& 14& 15\\
$d_n$ &
1& 1& 2& 1& 6& 2& 6& 3& 10&  2&  6&  2& 210& 30& 12\\
\hline
$n$ &
16& 17& 18&  19& 20&  21& 22& 23& 24&  25& 26& 27& 28& 29& 30\\ 
$d_n$ &
 3& 30& 10& 210& 42& 330& 30& 60& 30& 546& 42& 28&  2& 60&  4\\
\hline
\end{tabular}
\end{center}\end{table}
The first values of $d_n$ are given inTable~\ref{table:d_n}, for more values 
we refer to the Online Encyclopedia of Integer Sequences \cite{SloaneOEIS}, sequence
A338025.

\medskip
Our new algorithm is given in pseudo-code as Algorithm~\ref{alg:BCH} below. 
The following comments 
should provide sufficient evidence for the correctness of the algorithm. 
Throughout we use the notation $\coeff(v,X)$ for the coefficient of 
a word $v$ in an expression $X$.

\begin{algorithm}\label{alg:BCH} 
    \caption{Computation of BCH coefficients}
    \DontPrintSemicolon
    \KwIn{$(q_1,\dots,q_m)\in\mathbb{N}^m_{>0}$, $\mathit{Afirst}\in\{\mathrm{true}, \mathrm{false}\}$}
    \KwOut{Coefficient $h_w$ of $w=\AA^{q_1}\BB^{q_2}\cdots(\AA\lor\BB)^{q_m}$ (if $\mathit{Afirst}=\mathrm{true}$)
                        or $w=\BB^{q_1}\AA^{q_2}\cdots(\AA\lor\BB)^{q_m}$ (if $\mathit{Afirst}=\mathrm{false}$)
                        in the BCH series $\log(\ee^{\AA}\ee^{\BB})=\sum_w h_w w$}                        
    $n:=\sum_{i=1}^mq_i$\;
    $d:=n!\cdot d_n$\;
    $C:=(0)\in\mathbb{Z}^{n\times n}$\;
    $\mathit{Acurrent}:=\mathit{Afirst}$\;
    \If{$m$ {\rm is even}}{$\mathit{Acurrent}:={\rm\bf not}\,  \mathit{Afirst}$}
    $\ell:=0$\;
    \For{$i:=m,m-1,\dots,1$}
    { \label{line:9}
        \For{$r:=1,\dots,q_i$}
        {
            $\ell:=\ell+1$\;
            $h:=0$\;
            \uIf{$i=m$}
            {
                $h:=d/\ell!$\;
            }
            \ElseIf{$\mathit{Acurrent}$ {\rm\bf and} $i=m-1$}
            {
                $h:=d/(r!q_{m}!)$\;
            }
            $C_{1,\ell}:=h$\;
            \For{$k:=2,\dots,\ell-1$}
            {
                $h:=0$\;
                \For{$j:=1,\dots,r$}
                {
                    \If{$\ell>j$ {\rm\bf and} $C_{k-1,\ell-j}\neq 0$}
                    {
                        $h:=h+C_{k-1,\ell-j}/j!$\;
                    }
                }
                \If{$\mathit{Acurrent}$ {\rm\bf and} $i\leq m-1$}
                {
                    \For{$j:=1,\dots,q_{i+1}$}
                    {
                        \If{$\ell>r+j$ {\rm\bf and} $C_{k-1, \ell-r-j}\neq 0$}
                        {
                            $h:=h+C_{k-1,\ell-r-j}/(r!j!)$\;
                        }
                    }				
                }
                $C_{k,\ell}:=h$\;
            }
            $C_{\ell,\ell} := d$\;		
        }	
        $\mathit{Acurrent}:={\rm\bf not}\, \mathit{Acurrent}$\;
    }
    {\bf return} $(\sum_{k=1}^n(-1)^{k+1}C_{k,n}/k)/d$\;
    \end{algorithm}
\begin{description}
    \setlength{\itemsep}{2pt}
      \setlength{\parskip}{2pt}
    \item[\it Input:] 
    We consider the  word 
    $w=\AA^{q_1}\BB^{q_2}\cdots (\AA\lor\BB)^{q_m}$ 
    or
    $w=\BB^{q_1}\AA^{q_2}\cdots (\AA\lor\BB)^{q_m}$ 
    as a concatenation of $m$ alternating  blocks of  $\AA$s or $\BB$s whose lengths are 
    $q_1,\dots,q_m$.
    The boolean variable $\mathit{Afirst}$ indicates whether the first block
    is an $\AA$-block (or otherwise a  $\BB$-block).
    \item[\it Line 1:] $n=q_1+\ldots+q_m$ is the length of the word $w$.
    \item[\it Line 2:] $d=n!\,d_n$ is the smallest common denominator of all coefficients of degree $\leq n$
    of the BCH series.
    Here, $d_n$ can be computed using equation (\ref{eq:d_n}), or it can be taken from Table~\ref{table:d_n}.
    More generally, $d$ may be set to any nonzero multiple of $n!\,d_n$.
    \item[\it Line 3:] The array $(C_{k,\ell})\in\mathbb{Z}^{n\times n}$ 
    is initialized to zero. It will eventually contain 
     $C_{k,\ell} = d\cdot\coeff(v(\ell), Y^k)$, $k=1,\dots,\ell$, $\ell=1,\dots,n$,
    where 
    \begin{equation*}
    Y=\ee^{\AA}\ee^{\BB}-1=\sum_{i+j>0}\frac{1}{i!j!}\AA^i\BB^j,
    \end{equation*}
    and $v(\ell)=w_{n-\ell+1}\cdots w_n$ is the right subword of $w=w_1\cdots w_n$ 
    of length $\ell$
    starting
    at position $n-\ell+1$.
    
    \item[\it Lines 9--38:] 
    The outermost loop over $i$ processes the $m$ blocks in reverse order.
    The boolean variable $\mathit{Acurrent}$ indicates whether the  current $i$-th block
    is an $\AA$-block.
    
    \item[\it Lines 10--36:]The loop over $r=1,\dots,q_{i}$ combines
    with the outer loop over $i$ to form a loop over 
    $\ell=r+q_{i+1}+\ldots+q_m$ which processes the right subwords $v(\ell)$
    of lengths $\ell$.
    
    \item[\it Lines 12--18:]If $k=1$ then the current right subword $v(\ell)$ can only contribute 
    to $C_{k,\ell}=C_{1,\ell}=d\cdot\coeff(v(\ell),Y)$,
    if it has the form $v(\ell)=\AA^s\BB^t$ with $s+t=\ell$, and thus if it is contained in the last two blocks.
    %
    This contribution is $d/\ell!$ if
     $v(\ell)$ is entirely contained  in the last (i.e., the $m$-th) block
     such that $v(\ell)=\AA^\ell$ or $v(\ell)=\BB^\ell$,
     or it
    is $d/(r!q_m!)$ if $v(\ell)$ is contained in the last two blocks,
    where the next to last (i.e., the $(m-1)$-th) block has to be an
    $\AA$-block such that $v(\ell)=\AA^r\BB^{q_m}$.
    
    \item[\it Lines 19--34:]Let $u(\ell,j)$ denote the left subword of $v(\ell)$ of
    length $j$ such that
    $$v(\ell)=u(\ell,j)v(\ell-j),\quad j=0,\dots,\ell.$$
    For $k=2,\dots, \ell-1$ we have
    $$
    \coeff(v(\ell),Y\cdot Y^{k-1})
    =\sum_{j=0}^\ell\coeff(u(\ell,j),Y)\cdot\coeff(v(\ell-j),Y^{k-1}).
    $$
    Here we have $\coeff(v(\ell-j),Y^{k-1})=0$  for $j=\ell$.
    Similarly as before (cf.~lines 12--18), 
    we have
    $\coeff(u(\ell,j),Y)\neq 0$ only if $j\geq 1$ and if either $u(\ell,j)$ is 
    entirely contained in the current $i$-th block
    (or, more precisely, the current right subblock of length $r$ of the $i$-th block),
    or if it is entirely 
    contained in the  union of the $i$-th and the $(i+1)$-th block, where
    the   $i$-th block has to be an $\AA$-block.
    In the former case $u(\ell,j)=\AA^j$ or $u(\ell,j)=\BB^j$, $j=1,\dots, r$ such that 
    $\coeff(u(\ell,j),Y)=1/j!$, and in the latter case $u(\ell,j)=
    \AA^r\BB^{j_1}$ with $r+j_1=j$, $j_1=1,\dots,q_{i+1}$
    such that $\coeff(u(\ell,j),Y)=1/(r!j_1!)$.
    It follows 
    \begin{align*}
    C_{k-1,\ell} &= d\cdot \coeff(v(\ell),Y^k)
    = \sum_{j=1}^r\frac{1}{j!}C_{k-1,\ell-j}+
       f_i\sum_{j_1=1}^{q_{i+1}}\frac{1}{r!j_1!}C_{k-1,\ell-r-j_1} ,
    \end{align*}
    where $f_i=1$ if the $i$-th block is an $\AA$-block and 
    $f_i=0$ otherwise. This sum is computed in lines 20--33.
    Note that here $C_{k-1,\ell-j}$ and $C_{k-1,\ell-r-j_1}$ either are understood to
    be $=0$ if the second index is 0,
    or 
    they have already been computed
    during a previous pass of the loop over $\ell$ (i.e., the loops over
    $i$ and $r$ combined).
    
    Obviously the tests for $C_{k-1,\ell-j}\neq 0$ respectively
    $C_{k-1,\ell-r-j}\neq 0$ in lines 22 and 28 are  not strictly necessary, 
    but are there for efficiency reasons.
    
    \item[\it Line 35:] 
    A word of degree $\ell$ occurs in
    $Y^\ell$ if and only if for all $i=1,\dots, \ell$,
    its $i$-th letter corresponds to exactly one term of degree 1 
    of the $i$-th factor $Y$ of $Y^n$. Thus, $v(n)$ occurs in $Y^\ell$ exactly
    once and with coefficient $1$ so that $C_{\ell,\ell} = d\cdot\coeff(v(\ell),Y^\ell)=d$.
    
    Note that the case $k=\ell$ could also be handled by the above loop  over $k$. 
    Here it is handled separately for efficiency and because it is so simple.
    
    \item[\it Line 39:] The final result is computed 
    according to 
    $$\coeff(w,\log(\ee^{\AA}\ee^{\BB})) = 
    \sum_{k=1}^n\frac{(-1)^{k+1}}{k}\coeff(w, Y^k) 
    =\frac{1}{d}\sum_{k=1}^n\frac{(-1)^{k+1}}{k}C_{k,n}.
    $$
    \end{description}

    A key feature of the algorithm is that 
    it performs all of its computations in integer arithmetic. This
    means in particular, that the divisions
    in lines 14, 16, 23, 29, and the divisions by $k$ in line 39
    never have a remainder. (Of course, this does not apply to the
    final division by $d$ in line 39.) To prove this, it is not enough
    to know that the final result is a rational number with a denominator that
    is a divisor of $d=n!d_n$. It must also be ensured that during  the computation  no intermediate results not representable as integers can occur, 
    which cancel out at the end. 
    Without going into details, this holds because
    the computations of the algorithm follow  the 
    same pattern as
    the computation of the common denominator  $d=D_n=n!d_n$
     in the proof of  \cite[Proposition~1]{HHdenom},
     where the generic case is assumed and no cancellations are  taken into account.

\section{Performance of the algorithm}
\label{Sect:performance}
\begin{table}\label{table:performance}
    \caption{Running time of the calculation of all coefficients $h_w$ of 
    words $w=\AA^{q_1}\BB^{q_2}\cdots(\AA\lor\BB)^{q_m}$ 
    corresponding to all partitions $n=q_1+\dots+q_m$, $m\geq 1$,
    $q_1\geq q_2\geq\dots\geq q_m\geq 1$ of $n$ for all $n\leq N$.}
    \begin{center}
    \begin{tabular}{rrrrr}
    \hline
    $N$ & \#partitions &    data type & time (seconds)\\
    \hline
     19 &     2086 &  {\tt Int64} &          0.09 \\
     20 &     2713 & {\tt Int128} &          0.12 \\
     30 &    28628 & {\tt Int128} &          2.40 \\
     40 &   215307 & {\tt BigInt} &        609.72 \\
    \hline
    \end{tabular}
    \end{center}
\end{table}

To illustrate the high efficiency of Algorithm~\ref{alg:BCH} we use our
Julia implementation \cite{HHjulia} of the algorithm to compute a table of all coefficients
$h_w$ of the BCH series up to a given maximal degree $N$.
It is well known that the coefficients $h_w$ for 
$w=\AA^{q_1}\BB^{q_2}\cdots(\AA\lor\BB)^{q_m}$ are invariant under
permutations of the $q_i$ and invariant up to a sign
$(-1)^{q_1+\dots+q_m+1}$ under exchanging $\AA$ and $\BB$,
see \cite{G}. 
It thus suffices to calculate only the coefficients of the
words $w=\AA^{q_1}\BB^{q_2}\cdots(\AA\lor\BB)^{q_m}$ corresponding to all partitions $n=q_1+\dots+q_m$, $m\geq 1$,
$q_1\geq q_2\geq\dots\geq q_m\geq 1$ of $n$ for all $n\leq N$.
A table of coefficients corresponding to such partitions up to degree $N=20$ was first 
published in \cite{NewmanThompson}.

Table~\ref{table:performance} shows running times for calculating
BCH coefficients up to several maximal degrees $N$ on a
standard personal computer with a 3.0\,GHz Intel Core i5-2320 processor and 8\,GB of memory.
Here the timings are excluding  the time for  generating the partitions which can be neglected.
As mentioned in the introduction, to avoid integer overflow, appropriate 
integer data types have to be used, which are also shown in the table. 
(Here {\tt Int64} is the standard Julia 64-bit integer type, 
{\tt Int128} is the Julia 128-bit integer type, and
{\tt BigInt} is the Julia built-in arbitrary precision integer type.)
A comparison with the timings given in  \cite{WeyrauchScholz} 
shows that our implementation is
several orders of magnitudes faster than the implementations 
considered in \cite{WeyrauchScholz}.

\end{document}